\documentclass{article}
\usepackage{array, amsmath, amssymb, graphicx, stmaryrd,a4wide}
\usepackage[percent]{overpic}
\usepackage{subfig}
\usepackage{authblk}
\usepackage{amsfonts}
\usepackage{epsfig}
\usepackage{graphics} 
\usepackage{xspace}
\usepackage{bm}

\newcommand{\cpp}{{\mathfrak{p}}}

\newcommand{\tn}{|||}

\newcommand{\mcK}{\mathcal{K}}

\newcommand{\mcN}{\mathcal{N}}
\newcommand{\mcF}{\mathcal{F}}

\newcommand{\IR}{\mathbb{R}}
\newtheorem{rmk}{Remark}
\newtheorem{thm}{Theorem}

\newtheorem{prop}[thm]{Propositiom}

\begin{document}

\title{A cut finite element method with boundary value correction for the incompressible Stokes' equations}

% your contribution title if the original one is too long use an abbreviated title (for running head):
%\titlerunning{CutFEM boundary value correction for the Stokes' equations}

%\author{Erik Burman\inst{1} \and Peter Hansbo\inst{2} \and Mats G. Larson\inst{3}}
% Use \authorrunning{Short Title} for an abbreviated version of the author list (for running head):
\author[$\dagger$]{Erik Burman}
\author[$\ddag$]{Peter Hansbo}
\author[$\star$]{Mats G.\ Larson}
\affil[$\dagger$]{\small Department of Mathematics, University College London, London, UK--WC1E  6BT, United Kingdom}
\affil[$\ddag$]{\small Department of Mechanical Engineering, J\"onk\"oping University, SE-55111 J\"onk\"oping, Sweden}
\affil[$\star$]{\small Department of Mathematics and Mathematical Statistics, Ume{\aa} University, SE-901 87 Ume{\aa}, Sweden}
\maketitle

\begin{abstract}
We design a cut finite element method for the incompressible Stokes equations on curved domains.
The cut finite element method allows for the domain boundary to cut
through the elements of the computational mesh in a very general
fashion. To further facilitate the implementation we propose to use a piecewise
affine discrete domain even if the physical domain has curved boundary.
Dirichlet boundary conditions are imposed 
using Nitsche's method on the discrete boundary and the effect of the
curved physical boundary is accounted for using the boundary value
correction technique introduced for cut finite element methods in 
Burman, Hansbo, Larson, {\emph{A cut finite element method with boundary value correction}}, Math. Comp. 87(310):633--657, 2018.
\end{abstract}

\section{Introduction}
Let $\Omega$ be a domain in $\mathbb{R}^d$ with smooth boundary 
$\partial \Omega$ and exterior unit normal $n$. We will consider a cut
finite element method (CutFEM) for Stokes' problem on $\Omega$
with Dirichlet conditions. See \cite{BurClaHan15} and the references 
therein for an introduction to CutFEM. The Stokes problem takes the form: 
find $u:\Omega \rightarrow \IR^d$ and $p: \Omega \rightarrow \IR$
such that
\begin{alignat}{2}\label{eq:stokes}
-\Delta u + \nabla p&= f \qquad 
&& \text{in $\Omega$}
\\ \label{eq:stokesbc}
\nabla \cdot u & = 0 \qquad 
&& \text{in $\Omega$}
\\ \label{eq:divfree}
u &= g \qquad && \text{on $\partial\Omega$}
\end{alignat}
where $f\in [H^{-1}(\Omega)]^d$ and $g\in [H^{1/2}(\partial \Omega)]^d$ 
are given data. It follows from the Lax-Milgram Lemma that there 
exists a unique solution $u \in [H^1(\Omega)]^d$ and from Brezzi's
Theorem that  there 
exists a unique solution $p \in L^2_0(\Omega)$. We also have the 
elliptic regularity estimate
\begin{equation}\label{eq:ellipticregularity}
\|u\|_{H^{s+2}(\Omega)} + \|p\|_{H^{s+1}(\Omega)}\lesssim \|f\|_{H^s(\Omega)}, \qquad 
s \geq -1
\end{equation} 
Here and below we use the notation $\lesssim$ to denote less or 
equal up to a constant. The objective of the present paper is to
propose a cut finite element method for the problem
\eqref{eq:stokes}-\eqref{eq:divfree}. Unfitted finite element methods
for incompressible elasticity was first discussed in \cite{BBH09}, for
the coupling over an internal (unfitted) interface. The fictitious
domain problem for the Stokes' equations was then considered in
\cite{BH14,MLLR14,BCM15} and more recently the inf-sup stability for several different
well-known elements on unfitted meshes was proved \cite{GO16} and
further work on the Stokes' interface problem was presented in \cite{HLZ14}. The
upshot in the present contribution is that we, following \cite{BHL15},
use a piecewise affine representation of the physical boundary and
introduce a correction in the Nitsche formulation to correct for the
low order geometry error. This allows us to use for instance a
piecewise affine levelset for the geometry representation used in the
integration over the cut elements, while retaining the accuracy of a
(known) higher order representation of the boundary. This provides an
alternative to representing curved boundaries using isoparametric
mappings, see \cite{LPWL16} for an application of this technique to
the Stokes' equations. We will focus herein on the derivation of the CutFEM 
boundary value correction method for the Stokes' system, this is the topic 
of Section \ref{sec:method}. We then state some fundamental results 
(without proof) in Section \ref{sec:results} and finally we report some 
numerical examples in Section \ref{sec:numerics}.

\section{The cutFEM for Stokes' equations - derivation}\label{sec:method}
Here we will give the elements of the numerical modelling that leads
to the cut boundary value correction method for Stokes' equations. 

\subsection{The domain}
We let $\varrho$ 
be the signed distance function to $\partial \Omega$, negative on the inside and 
positive on the outside, and we let 
$U_\delta(\partial \Omega)$ 
be the tubular neighborhood $\{x\in \IR^d : |\varrho(x)| < \delta\}$ 
of $\partial \Omega$. Then there is a constant $\delta_0>0$ such 
that the closest point mapping $\cpp (x):U_{\delta_0}(\partial \Omega) 
\rightarrow \partial \Omega$ is well defined and we have the 
identity $\cpp (x) = x - \varrho(x)n(\cpp (x))$. We assume that $\delta_0$ is
chosen small enough that $\cpp(x)$ is well defined. 

\subsection{The mesh, discrete domains, and finite element spaces}
\begin{itemize}
\item Let $\Omega_0 \subset \IR^d$ be a convex polygonal domain such 
that $U_{\delta_0}(\Omega) \subset \Omega_0$, where 
$U_{\delta}(\Omega) := U_{\delta}(\partial \Omega) \cup 
\Omega$. Let $\mcK_{0,h}, h \in (0,h_0]$, 
be a family of quasiuniform partitions, with 
mesh parameter $h$, of $\Omega_0$ into shape 
regular triangles or tetrahedra $K$. We refer 
to $\mcK_{0,h}$ as the background mesh.

\item Given a subset $\omega$ of $\Omega_0$, let 
$\mcK_h(\omega)$ be the submesh defined by
\begin{equation}
\mcK_{h}(\omega) = \{K \in \mcK_{0,h} : \overline{K} \cap 
\overline{\omega} 
\neq \emptyset \}
\end{equation}
i.e., the submesh consisting of elements that intersect 
$\overline{\omega}$, and let 
\begin{equation}
\mcN_{h}(\omega) = \cup_{K \in \mcK_{h}(\omega)} K
\end{equation}
be the union of all elements in $\mcK_h(\omega)$. Below the $L^2$-norm 
of discrete functions frequently should be interpreted as the broken norm. For example for norms over $\mcN_{h}$ we have
\begin{equation}
\|v\|_{\mcN_{h}(\omega)}^2 := \sum_{K \in \mcK_{h}(\omega)} \|v\|_K^2
\end{equation}
\item Let $\Omega_h$, $h \in (0,h_0]$, be a 
family of polygonal domains approximating $\Omega$, possibly independent of
the computational mesh. We 
assume neither $\Omega_h \subset \Omega$ nor  $\Omega \subset
\Omega_h$, instead the accuracy with which $\Omega_h$ approximates
$\Omega$ will be important. 
\item Let the active mesh $\mcK_h$ be defined by
\begin{equation}
\mcK_{h} := \mcK_h(\Omega \cup \Omega_h)
\end{equation}
i.e., the submesh consisting of elements that intersect 
$\Omega_h\cup\Omega$, and let 
\begin{equation}
\mcN_{h} := \mcN_h(\Omega\cup \Omega_h)
\end{equation}
be the union of all elements in $\mcK_h$.
 \item Let $V^k_{0,h}$ be the space of piecewise continuous 
polynomials of order $k$ defined on $\mcK_{0,h}$ and let the 
finite element space $V^k_h$ be defined by
\begin{equation}
V^k_{h} :=\{ v_h: v_h := \tilde v_h\vert_{\mcN_h} \mbox{ for } \tilde
v_h \in  V^k_{0,h} \}
\end{equation}
\item To each $
\Omega_h$ we associate the normal $\nu_h:\partial \Omega_h \rightarrow
\mathbb{R}^d$, $|\nu_h|=1$, and the distance from $\partial \Omega_h$
to $\partial \Omega$, $\varrho_h:\partial \Omega_h \rightarrow \mathbb{R}$, such that if $\cpp_h(x,\varsigma):=x + \varsigma \nu_h(x)$ then
$\cpp_h(x,\varrho_h(x)) \in \partial \Omega$ for all $x \in \partial \Omega_h$. We will also assume that $\cpp_h(x,\varsigma)
\in U_{\delta_0}(\Omega)$ for all $x \in \partial \Omega_h$ and all
$\varsigma$ between $0$ and $\varrho_h(x)$. For conciseness we will 
drop the second argument of $\cpp_h$ below whenever it takes the value
$\varrho_h(x)$. 
We assume that the following assumptions are satisfied
\begin{equation}\label{eq:geomassum-a}
\delta_h := \| \varrho_h \|_{L^\infty(\partial \Omega_h)} = o(h^\zeta), 
\qquad  h \in (0,h_0]
\end{equation}
and 
\begin{equation}\label{eq:geomassum-c}
\| \nu_h - n\circ \cpp \|_{L^\infty(\partial \Omega_h)} = o(h^{\zeta-1}), 
\qquad  h \in (0,h_0]
\end{equation}
where $o(\cdot)$ denotes the little ordo and $\zeta \in \{1,2\}$. We also assume 
that $h_0$ is small enough for some additional geometric conditions to
be satisfied, for details see \cite[Section 2.3]{BHL15}.
\end{itemize}
\subsection{Numerical modelling} 
We now proceed to show how to obtain a boundary value correction
formulation for the Stokes' system.
\paragraph{Derivation.} Let $f=Ef$ and $u=Eu$ be the 
extensions of $f$ and $u$ from $\Omega$ to 
$U_{\delta_0}(\Omega)$. For $v \in V_h$ we 
have using Green's formula
\begin{align}
(f,v)_{\Omega_h} 
&= (f+\Delta u-\nabla p,v)_{\Omega_h} - (\Delta u - \nabla p,v)_{\Omega_h} 
\\ 
&=(f +\Delta u-\nabla p ,v)_{\Omega_h \setminus \Omega} 
+ (\nabla u,\nabla v)_{\Omega_h} - (p,\nabla v)_{\Omega_h}\\ \nonumber
&
\qquad - (\nu_h\cdot\nabla u + \nu_h p,v)_{\partial \Omega_h}
\end{align}
where we used the fact 
$f+ \Delta u-\nabla p = Ef - \Delta Eu - \nabla E p$, which is not in general equal 
to zero outside $\Omega$. Now the boundary condition $u=g$ on $\partial \Omega$ 
may be enforced weakly as follows
\begin{align}
(f,v)_{\Omega_h} 
&=(f +\Delta u-\nabla p ,v)_{\Omega_h \setminus \Omega} 
+ (\nabla u,\nabla v)_{\Omega_h} - (p,\nabla v)_{\Omega_h}\\ \nonumber
&\qquad - (\nu_h\cdot\nabla u + \nu_h p,v)_{\partial \Omega_h}
- (u\circ \cpp_h - g\circ \cpp_h,\nu_h\cdot \nabla v)_{\partial \Omega_h}
\\ \nonumber
&\qquad
+ \beta h^{-1} (u\circ \cpp_h - g\circ \cpp_h, v)_{\partial \Omega_h}  
\end{align}
Since we do not have access to $u\circ \cpp_h$ we use a Taylor 
approximation in the direction $\nu_h$
\begin{equation}\label{def:Taylor}
u\circ \cpp_h(x) \approx T_k(u)(x) 
:= \sum_{j=0}^k \frac{D_{\nu_h}^j u(x)}{j!}\varrho_h^j(x)
\end{equation}
where $D_{\nu_h}^j$ is the $j$:th partial derivative in 
the direction $\nu_h$.
Thus it follows that the solution to 
(\ref{eq:stokes})-(\ref{eq:divfree}) satisfies 
\begin{align}\label{eq:derivation-c}
&(f,v)_{\Omega_h} =(f +\Delta u-\nabla p ,v)_{\Omega_h \setminus \Omega} 
+ (\nabla u,\nabla v)_{\Omega_h} - (p,\nabla v)_{\Omega_h}
\\ \nonumber
&\qquad - (\nu_h\cdot\nabla u + \nu_h p,v)_{\partial \Omega_h}
\\ \nonumber
&\qquad - (T_k(u) - g\circ \cpp_h,\nu_h\cdot \nabla v)_{\partial \Omega_h}
+ \beta h^{-1} (T_k(u) - g\circ \cpp_h, v)_{\partial \Omega_h} 
\\ \nonumber
&\qquad  - (u\circ \cpp_h - T_k(u),\nu_h\cdot \nabla v)_{\partial \Omega_h}
+ \beta h^{-1} (u \circ \cpp_h - T_k(u), v)_{\partial \Omega_h} 
\end{align}
for all $v\in V_h$. Rearranging the terms we arrive at
\begin{align}
\nonumber
&(\nabla u,\nabla v)_{\Omega_h} - (p,\nabla v)_{\Omega_h}
-(\nu_h\cdot\nabla u + p \nu_h,v)_{\partial \Omega_h}
\\ \nonumber
&\qquad \qquad 
- (T_k(u),\nu_h\cdot \nabla v)_{\partial \Omega_h}
+ \beta h^{-1} (T_k(u), v)_{\partial \Omega_h} 
\\ \nonumber
& \qquad \qquad +(f + \Delta u-\nabla p, v)_{\Omega_h \setminus \Omega}
\\ \nonumber
&\qquad \qquad - (u\circ \cpp_h - T_k(u),\nu_h\cdot \nabla v)_{\partial \Omega_h}
+ \beta h^{-1} (u \circ \cpp_h - T_k(u), v)_{\partial \Omega_h} 
\\ \label{eq:derivation-d}
&\qquad=(f,v)_{\Omega_h}  
- (g\circ \cpp_h,\nu_h\cdot \nabla v)_{\partial \Omega_h}
+ \beta h^{-1}( g\circ \cpp_h, v)_{\partial \Omega_h} 
\end{align}
for all $v\in V^k_h$. The discrete method is obtained from this
formulation by dropping the consistency terms of highest order,
i.e. those on lines three and four of (\ref{eq:derivation-d}).
\paragraph{Bilinear Forms.}
We define the forms
\begin{align}\label{eq:a0}
a_0(v,w) &:= (\nabla v,\nabla w)_{\Omega_h} 
\\ \nonumber
&\qquad - (\nu_h\cdot \nabla v,w)_{\partial \Omega_h} 
- (T_k(v),\nu_h\cdot \nabla w)_{\partial \Omega_h} 
\\ \nonumber
&\qquad + \beta h^{-1}(T_k(v),w)_{\partial \Omega_h} 
\\ \label{eq:ah}
a_h(v,w)& := a_0(v,w) + j^k_-(v,w)
\\ \label{eq:bh}
b_\sigma (q,w)&:= (q,\nabla \cdot v)_{\Omega_h} - \sigma (q,v\cdot \nu_h)_{\partial \Omega_h}
\\ \label{eq:Sp}
s(y,q) &:= j^m_+(y,q) + \gamma_p\sum_{F \in \mcF_{h}}  h^{3} ([n_F \cdot
\nabla y],[n_F \cdot\nabla q])_F
\\ \label{eq:jh}
j^k_{\pm}(v,w) &:= \gamma_j\sum_{F \in \mcF_{h}} 
\sum_{l=1}^k h^{2l\pm1} ([D_{n_F}^l v],[D_{n_F}^l w])_F
\\ \label{eq:lh}
l_h(w) &:= (f,w)_{\Omega_h} 
- (g\circ \cpp_h ,\nu_h\cdot \nabla w)_{\partial \Omega_h} 
+ \beta h^{-1}(g\circ \cpp_h, w)_{\partial \Omega_h}                        
\end{align}
where $\gamma_j$, $\gamma_p$ and $\beta$ are positive constants. Here we used the notation:
\begin{itemize}
\item $\mcF_{h}$ is the set of all internal faces to 
elements $K \in \mcK_{h}$, i.e. faces that are not
  included in the boundary of the active mesh $\mcK_{h}$, that intersect the set  
$\Omega \setminus \Omega_h \cup \partial \Omega_{h}$, and 
$n_F$ is a fixed unit normal to $F\in \mcF_h$.
\item $D_{n_F}^l$ is the partial derivative of order 
$l$ in the direction of the normal $n_F$ to the face 
$F\in \mcF_{h}$.
\item $[v]|_F = v^+_F - v^-_F$, with $v_F^{\pm} 
= \lim_{s \rightarrow 0^+} 
v(x \mp s n_F)$, is the jump of a discontinuous function 
$v$ across a face $F\in \mcF_{h}$.
\item The stabilizing term $j_{h}(v,w)$ is introduced to extend the
  coercivity of $a_0(\cdot,\cdot)$ to all of $\mcN_h$ as we shall see
  below and similary for the pressure. Thanks to this property one may prove that the condition
  number is uniformly bounded independent of how $\Omega_h$ is
  oriented compared to the mesh following the ideas of
  \cite{Bu10,MLLR14}.
\item Observe the presence of the penalty coefficient
    $\beta$ in \eqref{eq:a0} and \eqref{eq:lh}. In order to guarantee 
    coercivity $\beta$ has to be chosen large enough and due to the 
    Taylor expansions we also have to require that $h \in (0,h_0]$ with 
    $h_0$ sufficiently small. 
\end{itemize}

\paragraph{The Method.} Find: $(u_h,p_h) \in W_h := [V^k_h]^d \times V^m_h$ such that 
\begin{equation}\label{eq:fem}
a_h(u_h, v) - b_1(p_h,v_h) + b_0(q_h,u_h) + s_p(p_h,q_h) = l_h(v),\qquad \forall (v,q) \in W_h
\end{equation}
where $a_h$ is defined in (\ref{eq:ah}), $b_0$ and $b_1$ in (\ref{eq:bh}) and $l_h$ in (\ref{eq:lh}). 

For the analysis below it will be convenient to use the compact
formulation:  find: $(u_h,p_h) \in W_h$  such that 
\begin{equation}
A_h[(u_h,p_h),(v,q)]+s(p_h,q) =  l_h(v),\qquad \forall (v,q) \in W_h
\end{equation}
where
\begin{equation}
A_h[(u_h,p_h),(v,q)]:=a_h(u_h, v) - b_1(p_h,v_h) + b_0(q_h,u_h)
\end{equation}
\begin{rmk}
Note that different forms $b_\cdot(\cdot, \cdot)$ are used in the moment and mass
equations and that
\begin{equation}
- b_1(p_h,u_h) + b_0(p_h,u_h) = (p_h, u_h \cdot \nu_h)_{\partial \Omega_h}
\end{equation}
It follows that the velocity pressure coupling term is skew-symmetric
only up to a boundary term that is essential for consistency. The
reason this term is omitted in the mass equation is that it is not
consistent and must either be improved using a special boundary value
correction, or omitted. For simplicity we here chose the latter option. 
\end{rmk}
\section{Theoretical results}\label{sec:results}
In this section we will report on some fundamental theoretical results
that hold for the formulation \eqref{eq:fem}.  Due to
space limitations the proof will be given elsewhere.

For the discussion we will introduce the following triple norms. We
will use the following norm defined for functions $(v,q)$ in $[H^2(\Omega)]^d
\times H^1(\Omega)$,
\begin{equation}
\tn (v,q) \tn_0 := \|v\|_{H^1(\Omega_0)}+ \|h^{-\frac12} v\|_{
\partial \Omega_h}+\|h \nabla p\|_{\Omega_0}+\|p\|_{\Omega_0}
\end{equation}
and an augmented version restricted to discrete
spaces or $(v,q)$ in $[H^2(\Omega)]^d
\times H^{\frac32+\epsilon}(\Omega)$, $\epsilon>0$,
\begin{align*}
\tn (v,q) \tn &:= \|v\|_{H^1(\Omega_0)}+ \|h^{-\frac12} v\|_{\partial \Omega_h} 
\\
&\qquad +\|h \nabla p\|_{\Omega_0}+\|p\|_{\Omega_0}+j_{-}(v,v)^{\frac12}+s(q,q)^{\frac12}
\end{align*}
\subsection{Inf-sup stability}
Key to discrete well-posedness and to the error analysis is the following inf-sup stability result that is
robust with respect to how the mesh intersects the interface. The main difficulty in the proof of
  this result is to handle the lack of skew symmetry between the terms
  $b_1$ and $b_0$. It follows however that the perturbation can be
  absorbed by the $L^2$-norm of the pressure and the boundary penalty
  term when $\beta$ is sufficiently large.
\begin{prop}\label{thm:infsup}
Let either $m=k$ and $\gamma_p>0$ or $k=2$ and $m=1$ and
$\gamma_p=0$ and assume that \eqref{eq:geomassum-a}-\eqref{eq:geomassum-c} hold with $\zeta=1$. Then there exists $\alpha >0$,  $h_0>0$ such that for
all $(v,q) \in W_h$, when $h<h_0$, there holds
\begin{equation}
\alpha \tn(v,q) \tn \leq \sup_{(w,y) \in W_h}
\frac{A_h[(v,q)(w,y)]+s(q,y)}{\tn (w,y) \tn}
\end{equation}
\end{prop}

\subsection{A priori error estimates}
In this section we will present an optimal error estimate in the norm
$\tn(\cdot,\cdot)\tn$. The proof of the estimate uses the classical structure of
inf-sup stability, Galerkin orthogonality, continuity of $A_h$,
estimation of geometry errors and finally approximability.
\begin{thm}
Let $(u,p) \in [H^{s}(\Omega)]^d \times H^{s-1}(\Omega)$, with
$s \ge 2$ be the solution to \eqref{eq:stokes} and assume that the
hypothesis of Proposition \ref{thm:infsup} are satisfied and that in addition \eqref{eq:geomassum-a}-\eqref{eq:geomassum-c} hold with $\zeta=2$. Let $(u_h,p_h) \in [V^k_h]^d \times V^m_h$ be
the solution of the finite dimensional problem \eqref{eq:fem}. Then
there holds
\begin{equation}
\tn(u - u_h,p-p_h)\tn_0 \lesssim  h^{\sigma} (|u|_{H^{\sigma+1}(\Omega)} + |p|_{H^{\sigma}(\Omega)})
\end{equation} 
where $\sigma = \min\{k,s-1\}$.
\end{thm}
\section{Numerical example}\label{sec:numerics}

In our numerical example we consider a two dimensional problem discretized by the lowest 
order (inf--sup stable) Taylor--Hood element: piecewise quadratic, continuous, approximation 
of the velocity and piecewise linear, continuous, approximation of the pressure, together with 
a piecewise linear approximation of the domain. We shall study the convergence with and without 
boundary modification.

We consider a problem from \cite{BH14} with exact solution (with $f=0$)
\begin{equation}
u_x =20x y^3, \quad u_y =5x^4-5y^4, \quad p=60x^2y-20y^3
\end{equation}
Our computational domain is a disc with center at the origin. The exact velocities 
are used as Dirichlet data on the boundary of the domain. Note that
since the exact solution is given everywhere, setting Dirichlet data
on the approximate boundary is not a problem in this (special) case;
to simulate the knowledge of data on the boundary only, we take the
boundary data from the edge of the exact domain and use as boundary
conditions on the approximate boundary, using the closest point
projection. We choose the method parameters $\gamma_1=\gamma_2=10^{-3}$, $\gamma_p=0$  
and $\beta=100$.

\begin{figure}[ht]
\begin{center}
\includegraphics[scale=0.15]{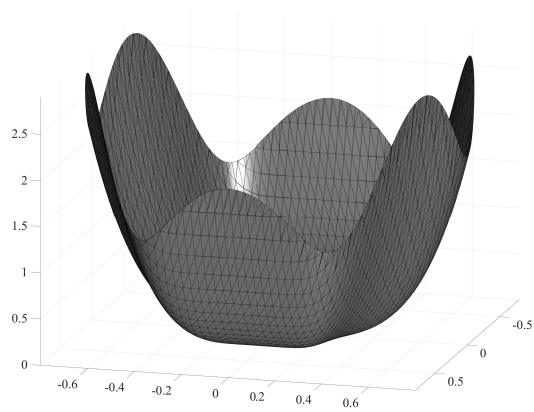}
\end{center}
\caption{Elevation of the norm of velocity (cut elements are triangulated for graphics purpose only).\label{fig:elev1}} 
\end{figure}

\begin{figure}[ht]
\begin{center}
\includegraphics[scale=0.25]{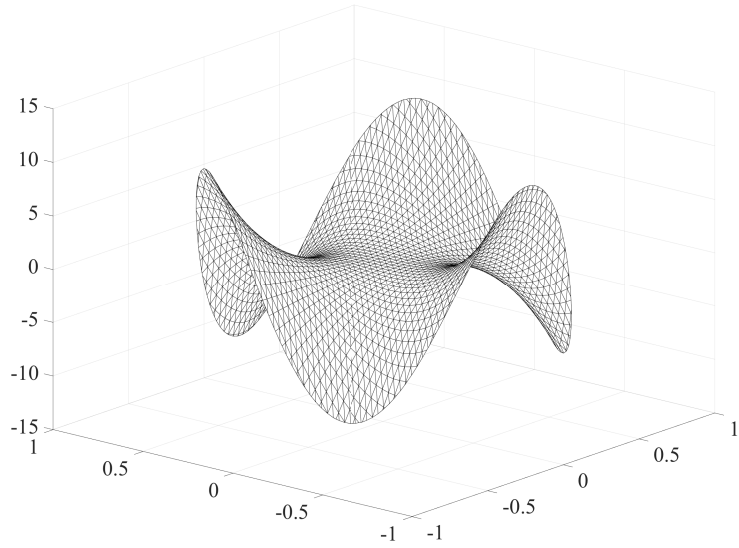}
\end{center}
\caption{Elevation of the pressure.\label{fig:elev2}} 
\end{figure}

\begin{figure}[ht]
\begin{center}
\includegraphics[scale=0.25]{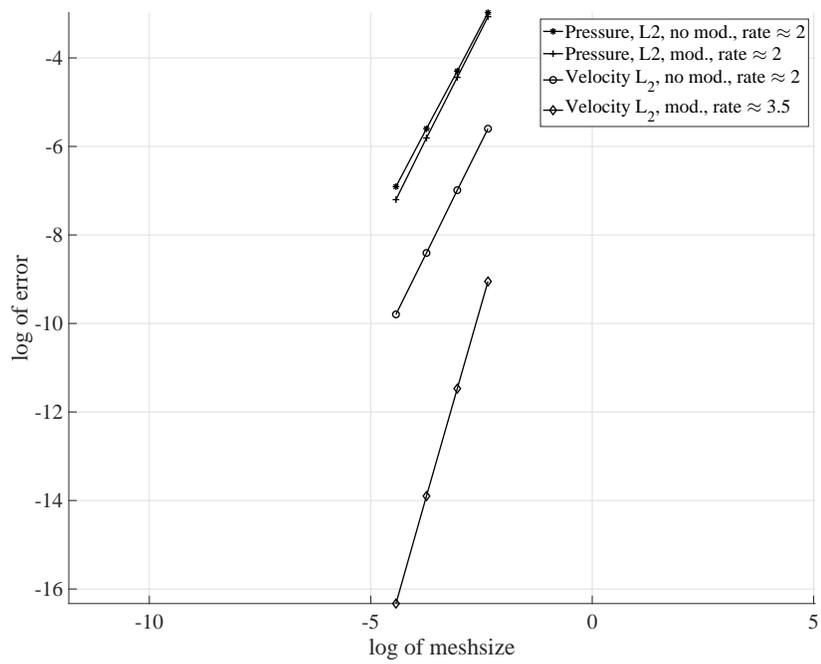}
\end{center}
\caption{Convergence results.\label{fig:conv}} 
\end{figure}

In Figure \ref{fig:elev1} and \ref{fig:elev2} we show elevations of the norm of velocity and the pressure, respectively. In Figure \ref{fig:conv} we show the convergence obtained with and without boundary modification. We note that without boundary modification we lose optimal convergence in 
velocities but retain optimal convergence for pressure, which is expected since the approximation 
of the boundary is piecewise linear leading to an $O(h^2)$ geometric consistency error. With boundary 
modification we recover optimal order convergence also for the velocity.

\subsection*{Acknowledgement}
This research was supported in part by EPSRC grant EP/P01576X/1, the Swedish Foundation for Strategic Research Grant No.\ AM13-0029, the Swedish Research Council Grants Nos.\  2013-4708, 2017-03911, and the Swedish Research Programme Essence.

 \vspace{5cm}

\end{document}